\documentstyle[12pt,twoside]{article}
\newtheorem{theorem}{Theorem}
\newtheorem{proposition}{Proposition}

\newtheorem{lemma}{Lemma}

\newtheorem{definition}{Definition}
\newtheorem{remark}{Remark}

\def\vsni{\vskip 0.2cm}

\def\per{{\rm Per}}

\def\N{I\!\!N}

\def\Z{I\!\!\!\!Z}

\def\t{\theta}

\def\S{\Sigma}
\def\S2{\Sigma_2}
\def\s{\sigma}

\def\ui{[0,1]}

\def\be{\begin{equation}}
\def\ee{\end{equation}}

\def\\{\hfill\break}
\def\={{ \; \equiv \; }}

\def\t_N{\tilde{\Z}_N}

\begin{document}
\date{}
%
%
%\begin{titlepage}
%
%
\title{On a set of numbers arising in the dynamics of unimodal maps}
\author{Stefano Isola \thanks{Dipartimento di Matematica e Informatica dell'Universit\`a
di Camerino and INFM, via Madonna delle Carceri, 62032 Camerino, Italy.
e-mail: $<$stefano.isola@unicam.it$>$.}}
\maketitle
\begin{abstract} 
\noindent
In this paper we initiate the study of the arithmetical properties of a set numbers which
encode the dynamics of unimodal maps in a universal way along with that of the corresponding topological zeta
function. Here we are concerned in particular with the Feigenbaum bifurcation. 
\end{abstract}
\vskip 1cm

\section{Preliminaries.}

We start by reviewing some basic ideas of (a version of) the {\sl kneading theory} for unimodal maps. For
related approaches and/or more details see
\cite{CE},
\cite{Dev},
\cite{deMvS}. 
\begin{definition}
A smooth map $f:\ui \to \ui$ is called {\sl unimodal} if it has exactly one critical point
$0<c_0<1$ and moreover $f(0)=f(1)=0$.
\end{definition}
For unimodal maps the orbit of the critical point 
$c_0$ determines in a sense the complexity of any other orbit. To be more precise, given
$x\in [0,1]$ we call  {\sl itinerary} of
$x$  with $f$ the sequence
$i(x) = s_1s_2s_3\dots$ where $s_{i} =0$ or $1$ according to
$f^{i-1}(x)< c_0$ or $f^{i-1}(x)\geq c_0$.
An important point is that such symbolic representation is in fact `faithful',
that is if
$s(x)=s(x')$ then
$x=x'$. Differently said, the partition of $[0,1]$ in the two semi-intervals
$P_0=[0,c_0)$ e $P_1=[c_0,1)$ is {\sl generating} for a unimodal map $f$ with critical point
$c_0$.
\vsni

\noindent
It is clear that if $s=i(x)$ is a sequence obtained as above then
$i(f(x))=\sigma (s)$ where
$\sigma$ denotes the left-shift: if $s=s_1s_2s_3\dots$ then $\sigma (s) =
s_2s_3s_4\dots $. The itinerary of the point $c_1=f(c_0)$ is called {\sl kneading sequence}
$K(f)$ of $f$. We say moreover that a given sequence $s$ of $0$ and $1$
is {\sl admissibile} for
$f$ if there is $x\in [0,1]$ such that
$i(x)=s$. A nice way to decide whether or not a given sequence is admissible 
amounts to establish an ordering on the itineraries which corresponds to ordering of the real
line. In this way, the admissible sequences are those which never become greater than the
kneading sequence when shifted.  To this end,
let us associate to a sequence $s=s_1s_2s_3\dots$ 
the number $\tau (s) \in [0,1]$ defined as
\be \label{tau}
\tau  = 0.t_1t_2t_3\dots =\sum_{k=1}^\infty {t_k\over 2^{k}},\qquad t_k = \sum_{i=1}^k s_i \,
({\rm mod}\, 2).
\ee
Equivalently, if we set
\be\label{eps1}
\epsilon_k =(-1)^{\sum_{i=1}^k s_i}
\ee
then $t_k$ and $\epsilon_k$ are related by
\be\label{taueps}
t_k ={1-\epsilon_k\over 2},\qquad \epsilon_k= 1-2\, t_k.
\ee
\begin{lemma}\label{tau} Given $x,y\in [0,1]$ we have
\begin{enumerate}
\item If $\tau (i(x)) < \tau (i(y))$ then $x<y$;
\item If $x<y$ then $\tau (i(x)) \leq \tau (i(y))$.
\end{enumerate}
\end{lemma}
\begin{remark}
The equality in 2) cannot be removed. Indeed, the existence of an attracting periodic orbit
typically implies the existence of an interval of points with the same itinerary.
On the other hand, a theorem due to Guckenheimer (see \cite{deMvS}) says that a unimodal map
$f$ has an attracting periodic orbit if and only if $K(f)$ is periodic. Viceversa, if $K(f)$ is
not periodic then implication  2) becomes: if $x<y$ then $\tau (i(x))< \tau (i(y))$.
This has important consequences. First of all: if $K(f)$ is not periodic and $K(f)
= K(g)$ then
$f$ and
$g$ are topologically conjugated.
\end{remark}
\vsni
\noindent
{\sl Proof of Lemma \ref{tau}.} Let us show the first part. Set
$i(x)=s_1s_2\dots$,
$i(y)=s'_1s'_2\dots$ and let $n=\min\{i\geq 1 : s_i\ne s'_i\}$ be the {\sl discrepancy} between
$i(x)$ and
$i(y)$. We proceed by induction in $n$. If $n=1$ the result is clear. Suppose it is true
for sequences with discrepancy $n-1$. We have $i(f(x)) = s_2s_3\dots $ and $i(f(y)) =
s'_2s'_3\dots$ Two cases are possbile:
either $s_1=0$ or $s_1=1$. If $s_1=0$ then $\tau(i(f(x))) < \tau (i(f(y)))$ because
applying
$f$ we don't modify the number of $1$'s before the discrepancy. Using the induction we then
have that
$f(x)<f(y)$. But since $f$ is increasing on $[0,c_0)$ we also have $x<y$. If $s_1=1$ then
$\tau (i(f(x))) >
\tau(i(f(y)))$
because there is a $1$ less among the symbols $s_2\dots s_n$. Therefore
by the induction we get
$f(x)>f(y)$ and since $f$ is decreasing on $(c_0,1]$ we see that $x<y$. The second
assertion follows similarly. $\Box$
\vsni
\noindent
An immediate consequence is the following  
\begin{theorem} \label{criterio} Every sequence $s$ such that
$$
\tau (\s (K(f))) \leq \tau (\s^m (s)) \leq \tau (K(f)),\qquad m\geq 0
$$
is admissibile and is the itinerary of a point in $[f(c_1),c_1]$.
\end{theorem}
In particular, 
$$
\tau (\s^m (K(f))) \leq \tau (K(f)),\qquad m\geq 0.
$$
A sequence $K$ with this property is said {\sl maximal}.
If moreover we consider a one-parameter family of unimodal maps $f_r$ so that $r\to f_r$ is continuous
on some real interval with respect to the ${\cal C}^1$ topology, then we can reformulate a theorem of 
Metropolis {\sl et al.}  (see \cite{deMvS}) by saying that every maximal sequence $K$
such that 
$$\tau (K({f_{r_a}}))\leq \tau (K) \leq \tau (K({f_{r_b}}))$$
is the kneading sequence of $f_r$ for some $r_a\leq r \leq r_b$. Notice that 
for $f_r([0,1]) \subseteq [0,1]$ one needs that $f_r(c_0)\leq 1$. In particular if $r=r_b$ then
$f_{r_b}(c_0)=1$ and $K({f_{r_b}})=1{\overline 0}$ (where ${\overline {s_1\dots s_l}}$ indicates
the unended repetition of the word ${s_1\dots s_l}$), to which it corresponds the number $\tau
(K({f_{r_b}}))=1$. At the other end-point we have $\tau (K({f_{r_a}}))=0$ (when $f^n(c_0)$ converges
monotonically to zero). We finally observe that given $q\in [0,1]$ we have
$$
\tau (q) =q \Longrightarrow \tau (\s (s)) = T(q)
$$
where $T:[0,1]\to [0,1]$ is the {\sl tent map} given by
$$
T(x) = \cases{ 2x &if $x< 1/2$, \cr 2(1-x) &if $x\geq  1/2$. \cr }
$$
Putting together these observations we obtain the following representation \cite{IP}: 

\vsni\noindent 

\begin{itemize}
\item The subset $\Lambda \subset [0,1]$ defined as
$$
\Lambda = \{ \tau \in [0,1] \; : \;  T^m(\tau) \leq \tau,\; \forall m\geq 0\}
$$
represents a {\sl universal encoding} for the dynamics of unimodal maps: those
having the same parameter $\tau$ have identical topological properties. In particular, every
$0$ in the binary expansion of $\tau \in \Lambda$ corresponds to a `forbidden word'
in the associated dynamics: let $K(f) = s_1s_2\dots $ and $\tau (K(f))=0.t_1t_2\dots$,
then if $t_{j}=0$ the word $s_1\dots {\hat s}_{j}$ (with ${\hat s}_{j}=1-s_j$) is
a forbidden word. 
Let $A=\{0,1\}$ be the alphabet and $A^*=\cup_{n\in \N}A^n$ the set of all possible finite words
written in the alphabet $A$. A word $u\in A^*$ of length $|u|=n$ is said {\sl
$f$-admissibile} if there is
$x\in [0,1]$ whose itinerary with $f$ up to the $n$-th letter coincides with $u$. The set
${\cal L}\subseteq A^*$ defined as
$$
{\cal L} = \{ u \in A^*\, ,Ê\, u \;\;\hbox{is $f$-admissibile} \}
$$
is the {\sl language} generated by $f$. The function
\be
p(n) =  \# \{ u\in {\cal L}, \,  |u|=n\}
\ee
is called the {\sl complexity function} of ${\cal L}$ and the limit
\be
h = \lim_{n\to \infty}{1\over n} \, \log \, p(n)
\ee
is the {\sl topological entropy}.
To summarize, the parameter $\tau$ furnishes a universal encoding in the sense that all unimodal maps
with  the same $\tau$ determine the same language ${\cal L} ={\cal L}(\tau)$ and, in particular, have the same
topological entropy
$h=h(\tau)$.
\end{itemize}
\begin{remark}
It is plain that the extremal situation in which  $T^m(\tau) = \tau$ for some $m> 0$ is that 
in which $\tau$ is a periodic point for the tent map $T$. In this case the kneading
sequence $K$ is periodic and so is the corresponding attractor.
This suggests that isolated points as well as `holes' in $\Lambda$ have to be related to periodic
attractors. In particular, there is a one-to-one correspondence between the holes in $\Lambda$  and the
periodic windows in the bifurcation diagram of unimodal maps, namely intervals in parameter space where the
topological entropy is constant
\cite{IP}.
\end{remark}
\subsection{Topological zeta function}
A great deal of information on the set of periodic points of a given map $f: [0,1]\to [0,1]$ can be stored into
the {\sl topological zeta function} of Artin and Mazur, defined as
\be
\zeta (f ,z) = \exp \sum_{n=1}^\infty {z^n\over n} \# \per_n (f).
\ee
For a unimodal map $f$ with parameter
$\tau$ the numbers $\# \per_n (f)$ are uniquely determined by the value of $\tau$ and we therefore write
$\zeta (\tau ,z) $. The series converges absolutely and uniformly for $|z|<e^{-h}$ and $z=e^{-h}$
is a singular point (e.g. a pole) of $\zeta (\tau,z)$.
This function can also be written as an Euler product noting that
$$
\sum_{n=1}^\infty {z^n\over n} \# \per_n (f) = \sum_{p=1}^\infty N(p)
\sum_{k=1}^\infty {z^{kp}\over k}=\log \prod_{p=1}^\infty (1-z^p)^{-N(p)},
$$
where $N(p)$ is the number of distinct periodic orbits of {\sl
prime period}
$p$. Hence we have
\be
\zeta (\tau,z) =\prod_{p=1}^\infty (1-z^p)^{-N(p)}.
\ee
The combinatorial features of the set of periodic orbits of a given map $f$ reflects onto the analytic
properties of 
$\zeta (\tau,z)$ in the complex plane.

\noindent
More specifically, it is not difficult to realize that from the work of Milnor and Thurston (\cite{MT}, Lemma
4.5 and Corollary 10.7) one can extract the following result
\begin{proposition}\label{zetatop} Let 
$ \Lambda \ni \tau=0.\,t_1t_2t_3\dots$. If the sequence $t_1t_2t_3\dots $ is eventually periodic or aperiodic
then 
\be
\zeta (\tau, z) = {1\over (1-z)\left(1+\sum_{k=1}^\infty \epsilon_k\, z^k\right)}.
\ee
If instead the kneading sequence is periodic and $\tau = 0.\, {\overline {t_1\dots t_{n}}}$
then
\be
\zeta (\tau, z) = {1\over (1-z)\left(1+\sum_{k=1}^{n-1} \epsilon_k\, z^k\right)},
\ee
where the numbers $\epsilon_k$ are defined in (\ref{eps1})-(\ref{taueps}).
\end{proposition}
\vsni
\noindent
We now list some examples in which the zeta function can be written in closed form by means of Lemma 
\ref{zetatop}.
\begin{itemize}
\item The number $\tau =1 = 0.1{\overline {1}}$ corresponds to the situation where the critical point gets
mapped to the origin in two steps, and yields
$$
\zeta (1, z) = {1\over 1-2z}, \qquad h(1) = \log 2.
$$
\item the number $\tau =5/6 = 0.1{\overline {10}}$  corresponds to the situation where the critical point gets
mapped to the fixed point in three steps ({\sl band merging}). In this case we find
$$
\zeta ({5/ 6}, z) = {1+z\over (1-z)(1-2z^2)},\quad h({5/ 6}) = \log \sqrt{2}.
$$
\item the number $\tau =6/7 = 0.{\overline {110}}$ corrisponds to the opening of the period three window.
The last orbit in the Sarkovskii order settles down and thus there are periodic orbits of any period. Here we
get
$$
\zeta ({6/7}, z) = {1\over (1-z)(1-z-z^2)},\quad h({6/7}) = \log{\sqrt{5}-1\over 2}.
$$
\end{itemize}
In the examples above the number $\tau$ was always rational. In the next Section we show a situation
leading to a trascendental irrational $\tau$.
A systematic study of the arithmetical properties of the numbers in $\Lambda$, along with their relation with
the dynamics, is far from being reached. In particular, the question of what is the most irrational $\tau$
(and to which chaotic state it corresponds) is open. In Section 2 we shall study the above
quantities for the Feigenbaum bifurcation but in order to get a self-contained exposition we first recall some
standard notions (for details see
\cite{Dev}).
\subsection{Kneading theory and renormalization}
\vsni\noindent
Let $f:[0,1]\to [0,1]$ be a unimodal map with a unique fixed point
$b$ in the interval
$(c_0,1)$, so that $f'(b) <0$.
Let $a$ the (unique) point 
in $(0,c_0)$ such that $f(a)=b$ and set
$J=[a,b]$. 
Consider the linear map $L$
defined by 
\be
L (x) = {1\over a - b}(x-b).
\ee
It expands $J$ to $\ui$ reversing its orientation. The inverse map is 
\be
L^{-1} (x) = (a - b)x+b.
\ee
The renormalization operator ${\cal R}$ is thus defined as 
\be\label{rino}
({\cal R} f)(x) = L\circ f^2 \circ L^{-1} (x).
\ee
Plainly $({\cal R} f)(0)=({\cal R} f)(1)=0$ and
$c_0$  is the only critical point of ${\cal R} f$. Moreover, $2$-periodic points for
$f$ become fixed points of ${\cal R} f$.

\noindent
Now let
$K(f) = s_1s_2s_3\dots$ be the kneading sequence of $f$. The following properties are easily verified
(see \cite{Dev}):
\begin{enumerate}
\item if ${\cal R}f$ is defined and unimodal then
$s_{2k+1}=1$,
$\forall k\geq 0$;
\item $K({{\cal R}f})={\hat s}_{2}{\hat s}_{4}{\hat s}_{6}\dots $. In other words, one can define a
renormalization operator on sequences acting as (with slight abuse we keep using the same symbol):
\be
{\cal R} (s_1s_2s_3\dots ) = {\hat s}_{2}{\hat s}_{4}{\hat s}_{6}\dots
\ee
\item if both ${\cal R}f$ and ${\cal R}^2f$ are unimodal then $s_{4k+2}=0$;
\item if ${\cal R}^lf$ is unimodal for $l\leq n$ then all symbols $s_j$ with $j=2^nk+2^{n-1}$ are determined;
\item since the numbers $2^nk+2^{n-1}$ exhaust all even numbers as $n$ varies in $\N$ 
it follows that if ${\cal R}^nf$
is unimodal for each $n\geq 1$ then all symbols of $K(f)$ are determined.
\end{enumerate}
How $K(f)$ looks like for an infinitely renormalizable unimodal map, that is a map $f$ such that ${\cal R}
f=f$?  
\vsni\noindent
Set 
\begin{eqnarray}
K_1 &=& {\overline {1}}\nonumber \\
K_2 &=& {\overline {10}}\nonumber \\
K_3 &=& {\overline {1011}}\nonumber \\
K_4 &=& {\overline {10111010}}\nonumber \\
K_5 &=& {\overline {1011101010111011}}\nonumber
\end{eqnarray}
and more generally
$K_{j+1}$ is obtained from $K_j$ by applying one of the following equivalent procedures: 
\begin{itemize} 
\item
duplicating the repeating sequence and reversing the last symbol; 
\item doubling all indices, reversing the resulting symbols (all with even index) and
inserting a
$1$ at each position with odd index; 
\item applying the {\sl Feigenbaum substitution} $1\to 10$ and $0\to 11$ (the symbol $1$ being the {\sl
prefix}) to the repeating sequence.
\end{itemize}
By construction $K_j$ has period $2^j$  with an odd number of $1$'s.
We also have that
\be
{\cal R}(K_{j+1})=K_j,\qquad j\geq 1.
\ee
Therefore the limit sequence
\be
K_\infty = \lim_{j\to \infty} K_j\, = \, 1011\, 1010\, 1011\, 1011\, \, 1011\, 1010\, 1011\, 1010\,\dots
\ee
is aperiodic and invariant under renormalization (the latter can be interpreted as a self-similarity
property):
\be
{\cal R}(K_{\infty})=K_\infty.
\ee
For all $j\geq 0$, $K_j$ is a prefix of $K_\infty$.
Finally, one easily verifies that $K_j$ is the kneading sequence of a unimodal map
having a periodic attractor of period $2^j$, whereas
$K_\infty$
is that of an infinitely renormalizable map.
\vsni
\noindent
Inspection of the sequences $K_n$ suggests that the asymptotic frequencies of the symbols 
$0$ and $1$ appearing in $K_\infty$ are $1/3$ and $2/3$ respectively.
To check this, we shall use a standard technique in the theory of substitution (see \cite{PF}) : let
$\phi$ be the substitution $\phi (1) =10$ and
$\phi (0)=11$ considered above and $N_i(\phi(j))$ be the number of occurrences of the symbol $i=0,1$
in the word
$\phi (j)$. The asymptotic frequency of $i$ in $K_\infty$ is then given by
\be
f_i = \lim_{n\to \infty} {N_i(\phi^n(1))\over 2^n}, \qquad i=0,1,
\ee
where we have used the fact that $|K_n|=2^n$. To compute $f_i$ we construct the matrix
\be
M=[N_i(\phi(j))]_{i,j\in \{0,1\}}.
\ee
A short reflection yields
\be
M^n=[N_i(\phi^n(j))]_{i,j=0,1},
\ee
and thus, setting $u=(0,1)$, we get
\be\label{fre}
f_i =  \lim_{n\to \infty} {(M^nu)_{i}\over 2^n}\cdot
\ee
From Perron-Frobenius theorem we have that $M$ has a simple positive eigenvalue of maximal modulus
$\lambda$ to which it corresponds an eigenvector with strictly positive components.
In our case we find
\be
M=\pmatrix{0  &1  \cr
            2  &1    \cr}
\ee
whose eigenvalues are $2$ and $-1$.
The normalized eigenvector corresponding to the leading eigenvalue is $v=(1/3,2/3)$. From (\ref{fre}) one
deduces that
$f_i=v_i$, $i=0,1$, which are the claimed frequencies.
\vsni
\noindent
\begin{remark}
One may consider the sequence
$K_\infty$ as an element of $\{0,1\}^{\N}$ and observe that the continuous injective map
$T:
\{0,1\}^{\N}\to 
\{0,1\}^{\N}$ defined as follows: if $\omega =111\dots$ then $T\omega = 000\dots$; if $\omega
=1\dots 10\dots$ then $T\omega = 0\dots 01\dots$; if $\omega
=0\dots$ then $T\omega = 1\dots$, acts a (right) translation on $K_\infty$. Therefore $T$ leaves invariant the
space $X={\overline {\{T^jK_\infty\}_{j\geq 0}}}$. The map $T:X\to X$ is called {\rm dyadic adding machine}.
\end{remark}

\section{Arithmetics of the Feigenbaum bifurcation}
\vsni
\noindent
We now look at the values of the parameter $\tau$ corresponding to the kneading sequences arising in the
period doubling scenario discussed in the preceeding Section. 

\noindent
Set
$\tau_j =
\tau (K_j)$. We find
\begin{eqnarray}
\tau _1 &=& 0.{\overline {10}}\nonumber \\
\tau _2 &=& 0.{\overline {1100}}\nonumber \\
\tau _3 &=& 0.{\overline {11010010}}\nonumber \\
\tau _4 &=& 0.{\overline {1101001100101100}} \nonumber\\
\tau _5 &=& 0.{\overline {11010011001011010010110011010010}} \nonumber 
\end{eqnarray}
and $\tau_{j+1}$ is obtained from $\tau_j$ by applying the rule
\be\label{rule}
\tau _j = 0.\, {\overline {t_1\dots t_{2^j}}}\; \Longrightarrow \; 
\tau _{j+1} = 0.\, {\overline {t_1\dots t_{2^j-1}{\hat t_{2^j}} \,{\hat t_{1}}  \dots {\hat t_{2^j-1}}
t_{2^j}}},
\ee
or, alternatively, by the following substitution: let 
$$
a=00,\quad b=01,\quad c=10,\quad d=11,
$$
then
\be\label{sub}
a \to ac,\quad b \to ad,\quad c \to da,\quad d \to db.
\ee
It is easy to check that $\tau_j \in \Lambda$, $\forall j\geq 1$.  They form an
increasing sequence: 
$$
\tau_1 < \tau_2 < \tau_3 < \cdots
$$
and satisfy
\be
{\tilde {\cal R}}(\tau_{j+1})=\tau_j\quad \hbox{where}\quad {\tilde {\cal R}}(0.t_1t_2t_3\dots
):=0.t_2t_4t_6\dots 
\ee
For each $j\geq 1$, $\tau_j = 0.\, {\overline {t_1\dots t_{2^j}}}$ is the rational
number given by
\be\label{rational}
\tau_j = {2^{2^{j}}\over 2^{2^{j}}-1}\; \sum_{k=1}^{2^{j}} {t_k\over 2^k}
\ee
We have
$$
\tau_1={2\over 3},\quad \tau_2={4\over 5},\quad \tau_3={14\over 17},\quad
\tau_4={212\over 257}, \quad \tau_5={54062\over 65537},
$$
$$
\tau_6={3542953172\over 4294967297},\qquad \tau_7={15216868001456509742\over 18446744073709551617}
$$
By (\ref{rule}) and (\ref{rational}) the following recursive law is in force:
\be
\tau_j={p_j\over q_j} \Longrightarrow \tau_{j+1}={p_{j+1}\over q_{j+1}}={2+p_j(q_j-2)\over
2+q_j(q_j-2)}
\ee
where all fractions are in lowest terms.
From the above it follows  $q_{j+1}-1=(q_j-1)^2$ and thus $q_j=2^{2^{j-1}}+1$. Note that the above
recursion can be written in the form
\be
p_1=2, q_1=3,\quad p_{j+1}=2+(2^{2^{j-1}}-1)p_j,\quad q_{j+1}=2+(2^{2^{j-1}}-1)q_j,\quad j\geq 1
\ee
This yields
\be
q_j-p_j= (2^{2^{j-2}}-1)(q_{j-1}-p_{j-1})= \cdots = \prod_{k=0}^{j-2}(2^{2^{k}}-1)
\ee
and recalling that $q_j=2^{2^{j-1}}+1$ we get
$p_j= 2^{2^{j-1}}+1- \prod_{k=0}^{j-2}(2^{2^{k}}-1)$.
We thus find the expression
\be
\tau_j = 1-{ \prod_{k=0}^{j-2}(2^{2^{k}}-1) \over 2^{2^{j-1}}+1}=1-{ \prod_{k=0}^{j-1}(1-2^{-2^{k}}) \over
2(1-2^{-2^{j}})}
\ee
and 
\be
\tau_{j+1}-\tau_j  =\left({2\over2^{2^j}+1}\right)
\tau_j.
\ee

\noindent
The number
\be\label{sviluppo2}
\tau_\infty =\lim_{j\to \infty} \tau_j =1-{1\over 2}\prod_{k=0}^{\infty}(1-2^{-2^{k}}) =0. 11010011\, 00101101
\, 00101100  \,11010011 \dots
\ee
satisfies $\tau_\infty = \tau (K_\infty)$ and is plainly irrational
(since $K_\infty$ is aperiodic). One easily recognizes the Thue-Morse sequence beginning in $0$, that is the
fixed point of the substitution $0 \to 01$ and $1\to 10$ with prefix $0$\footnote{By the way, we have shown
the following result: \begin{proposition} Let $\xi : \{0,1\}^{\N} \to \{0,1\}^{\N}$ be the map defined as $(\xi
s)_k=\sum_{i=1}^k s_i \,
({\rm mod}\, 2)$. Let $u$ be the fixed point of the Feigenbaum substitution $1\to 10$ and $0\to 11$ with
prefix $1$ and $w$ be the fixed point of the Thue-Morse substitution $0 \to 01$ and $1\to 10$ with prefix $0$.
Then $0\, \xi (u) = \xi (0u)= w$.
\end{proposition}}. 
It enjoys the invariance property
\be
{\tilde {\cal R}}(\tau_{\infty})=\tau_\infty
\ee
 which can also be expressed in the form
\be\label{invar}
\tau_\infty = \sum_{k=1}^{\infty} {t_k\over 2^k} = \sum_{k=1}^{\infty} {t_{2^lk}\over 2^k},
\qquad \forall l\geq 0.
\ee
Thus, for instance, $t_k=1$ whenever $k=2^{\ell}$ for some $\ell \geq 0$. More specifically,
we have
\begin{proposition}
For an integer $p\geq 1$ set
$$
s(p) = \sum_{i\geq 0} n_i \,({\rm mod}\, 2) \quad \hbox{if}\quad p = \sum_{i\geq 0} n_i \,2^i,\quad
n_i\in\{0,1\}.
$$
Let $\tau_\infty=0.t_1t_2\dots$.
Then
$t_k = s(p)$ whenever $k=p\cdot 2^{\ell} $ for some $\ell \geq 0$ and $p\geq 1$ odd.
\end{proposition}
\noindent
{\sl Proof.} Due to (\ref{invar}) it will suffice to show by induction over $r$ the following property:
$P_r = \{ p\; {\rm odd}\; {\rm  and} \; p\leq 2^r \Rightarrow t_k=s(p)\}$. Note that $P_0$ is obvious.
Consider an odd $p'$ such that $2^r<p'\leq 2^{r+1}$. Then $p'=2^r+p$ with $1\leq p\leq 2^r$ and $p$ odd. 
Then $s(p')=s(p)+1\, ({\rm mod}\, 2)$ and by the above $P_r \Rightarrow P_{r+1}$. $\Box$

\vsni
\noindent
Furthermore, from (\ref{sub}) we see at once that the symbols
$0$ and
$1$ both appear in
$\tau_\infty$ with frequency $1/2$. One may wonder if $\tau_\infty$ is a {\sl normal} number, in the sense of
Borel. That means that in its dyadic expansion (\ref{sviluppo2}) the asymptotic frequency
of any word of length $n$ is 
$2^{-n}$. On the other hand, reasoning as for the sequence $K_\infty$ (and using the substitutions
(\ref{sub})) it is not difficult to verify that the frequency of the pairs $00$, $01$, $10$ and
$11$ are
${1\over 3}$, 
${1\over 6}$, 
${1\over 6}$ and
${1\over 3}$. Therefore  $\tau_\infty$ is not a normal number. In fact $\tau_\infty$ is trascendental, as is
shown by Mahler in
\cite{Ma} (see also \cite{Dek}, \cite{AZ}, \cite{FM}) .

\vsni
\noindent
We end this digression with some partial insigth into the structure of the continued fraction
expansion of $\tau_\infty$. 

\noindent
Recall that any number $\tau \in [0,1]$ can be expanded as 
\be\label{cfe}
\tau = {1\over \displaystyle a_1 + {1\over \displaystyle a_2 + {1\over\displaystyle a_3
+\cdots }}}\equiv [a_1,a_2,a_3,\dots]
\ee 
where the $a_i$'s are integers. Successive truncations of this expansion yields a sequence of rational
numbers
\be\label{approx}
{r_n  \over s_n }=[a_1,a_2,a_3,\dots, a_n] 
\ee 
which are called {\sl convergents} of $\tau$ (see \cite{Kh}).

\noindent
Now, the problem we are interested in is the following: are the continued fraction
expansions of the numbers $\tau_j$ predictable (i.e. have a definite pattern) as their binary expansions do?
The expansions of the
first eight
$\tau_j$'s are
\begin{eqnarray}
\tau _1 &=& [1,2]\nonumber \\
\tau _2 &=& [1,4]\nonumber \\
\tau _3 &=& [1,4,1,2]\nonumber \\
\tau _4 &=& [1,4,1,2,2,6] \nonumber\\
\tau _5 &=& [1,4,1,2,2,6,2,1,2,9,1,2] \nonumber\\
\tau _6 &=& [1, 4, 1, 2, 2, 6, 2, 1, 2, 9, 1, 2, 2, 1, 1, 21, 1, 10, 2, 1, 1, 1, 5]\nonumber\\
\tau _7 &=& [1, 4, 1, 2, 2, 6, 2, 1, 2, 9, 1, 2, 2, 1, 1, 21, 1, 10, 2, 1, 1, 1, 4, 1, 2, 
29, 1, 24, 1,\nonumber\\
&& 1, 7, 11, 3, 2, 5, 1, 1,
1, 89]\nonumber\\
\tau _8 &=& [1, 4, 1, 2, 2, 6, 2, 1, 2, 9, 1, 2, 2, 1, 1, 21, 1, 10, 2, 1, 1, 1, 4, 1, 2, 
29, 1, 24, 1,\nonumber \\
&& 
1, 7, 11, 3, 2, 5, 1, 1, 1,
88, 1, 1, 1, 6, 1, 1, 33, 2, 6, 1, 24, 1, 5, 212, 2, 1, 10, 1,\nonumber \\ &&   3, 11, 2, 1, 2, 1,
10, 1, 1, 2, 3, 2549, 1, 2]
\nonumber
\end{eqnarray}
A direct inspection suggests that there is a subsequence $n_j$ of the integers so that if
$\tau_j=[a_1,\dots ,a_{n_j}]$ then 
\be
\tau_{j+1}=
\cases{[a_1,\dots ,a_{n_{j}}-1,b_{n_{j}+1},\dots ,b_{n_{j+1}}] &if $n_j$ is odd,\cr 
[a_1,\dots ,a_{n_{j}}, \;\;\;\;\;\;b_{n_{j}+1},\dots ,b_{n_{j+1}}]&if $n_j$ is even, \cr }
\ee
for some  $b_{n_{j}+1},\dots ,b_{n_{j+1}}$. The sequence $n_j$
for $1\leq j\leq 12$ is
$$
2,\; 2,\; 4,\; 6,\; 12,\;  23,\; 39,\; 71,\; 121,\; 253,\; 528,\;  1129
$$

\noindent
Unfortunately we are not able to say much more. In particular it is not clear what kind of relation could
be established between the $\tau_j$'s and the
convergents of
$\tau_\infty$.
Note that Shallit obtained in
\cite{S} a rather complete description of the patterns arising for irrational numbers of the type
$\sum_{k\geq 0}u^{-2^k}$, $u$ an integer. On the other hand, a high-temperature-like expansion of
the product appearing in (\ref{sviluppo2}) yields the expression
\be
\tau_\infty = 1- {1\over 2}\left[ 1-\sum_{\ell =1}^\infty {(-1)^\ell \over \ell !}
\sum_{k_1\ne k_2\ne \cdots \ne
k_\ell, \atop k_i\geq 0} 2^{-\sum_{i=1}^\ell 2^{k_i}}\right]
\ee 
of which the numbers studied by Shallit are just the first order $(\ell =1)$ term with $u=2$. We conclude with
a brief description of the topological zeta functions arising in this situation.

\vsni
\noindent
{\sc Zeta functions.}
For the values $\tau=\tau_j$ considered above, we get the polynomial zeta function
\be
1/\zeta (\tau_j, z) = (1-z) \prod_{n=0}^{j} (1-z^{2^n}),
\ee
whose zeroes are all on the unit circle $|z|=1$. Moreover we have
$\zeta (\tau_j, z) \to \zeta (\tau_\infty, z)$ when $j\to \infty$, where
\be
1/\zeta (\tau_\infty, z) = (1-z) \prod_{n=0}^{\infty} (1-z^{2^n}).
\ee
From Sarkovskii theorem (see \cite{BGMY}) it follows that 
$h(\tau_j)=0$ for all $j\geq 0$. Also 
$h(\tau_\infty)=0$ (but for any $\tau \in \Lambda$ with $\tau > \tau_\infty$ we have $h(\tau)>0$).
Put 
\be\label{produ}
\Xi (z) = \prod_{n=0}^{\infty} (1-z^{2^n}).
\ee
This function satisfies the functional equation 
\be
\Xi (z)=(1-z)\,\Xi (z^2)
\ee
from which we see that if $\Xi (z)=0$ then $|z|=1$. In particular, given $m\geq 1$ and $k=0,1,\dots , 2^l -1$ 
all factors of the product defining $\Xi (z)$ corresponding to  $n\geq m$
vanish at
$z=e^{2\pi ik/ 2^l}$. Therefore the zeroes of $\Xi (z)$ are dense on the unit circle. 
We then have that the radius of convergence of  $\zeta (\tau_\infty, z)$ is equal to
$1$ and that the unit circle is a (opaque) natural boundary for this function.

\noindent
Finally, expanding the product (\ref{produ}) we
get
\begin{eqnarray}
\Xi (z) &=& 1 - z - z^2 + z^3 - z^4 + z^5 + z^6 - z^7 - z^8 + z^9 +\nonumber \\
&&\qquad \qquad + z^{10} - z^{11} + z^{12} - z^{13} - z^{14} + z^{15} - z^{16} + \cdots \nonumber
\end{eqnarray}
from which we see that if
$\tau_\infty = 0.\,t_1t_2t_3\dots $ then the coefficient of $z^k$ with $k\geq 1$ in the above expansion
is but the number $\epsilon_k= 1-2\, t_k$ defined in (\ref{taueps}), in agreement with Proposition
\ref{zetatop}. In turn, we notice that the number $\tau_\infty$ can be written as 
\be
\tau_\infty =1-{1\over 2}\, \Xi \left( {1\over 2}\right).
\ee

\end{document}